# Optimum Weight Selection Based LQR Formulation for the Design of Fractional Order $PI^\lambda D^\mu$ Controllers to Handle a Class of Fractional Order Systems

Saptarshi Das[a,b], Indranil Pan[a]
a) Department of Power Engineering, Jadavpur University, Salt Lake Campus, LB-8, Sector 3, Kolkata-700098, India.
b) Communications, Signal Processing and Control Group, School of Electronics and Computer Science, University of Southampton, Southampton SO17 1BJ, United Kingdom.
Email: saptarshi@pe.jusl.ac.in, s.das@soton.ac.uk

Kaushik Halder[a,c], Shantanu Das[d] and Amitava Gupta[a]
c) Department of Electronics and Instrumentation Engineering, National Institute of Science and Technology, Palur Hills, Berhampur-761008, Orissa, India.
d) Reactor Control Division, Bhabha Atomic Research Centre, Mumbai-400085, India.

*Abstract*—A weighted summation of Integral of Time Multiplied Absolute Error (ITAE) and Integral of Squared Controller Output (ISCO) minimization based time domain optimal tuning of fractional-order (FO) PID or $PI^\lambda D^\mu$ controller is proposed in this paper with a Linear Quadratic Regulator (LQR) based technique that minimizes the change in trajectories of the state variables and the control signal. A class of fractional order systems having single non-integer order element which show highly sluggish and oscillatory open loop responses have been tuned with an LQR based FOPID controller. The proposed controller design methodology is compared with the existing time domain optimal tuning techniques with respect to change in the trajectory of state variables, tracking performance for change in set-point, magnitude of control signal and also the capability of load disturbance suppression. A real coded genetic algorithm (GA) has been used for the optimal choice of weighting matrices while designing the quadratic regulator by minimizing the time domain integral performance index. Credible simulation studies have been presented to justify the proposition.

## I. Introduction

Modern optimal control theory proposes several analytical tools to design not only control strategies satisfying desirable characteristics according to the designer's specifications, but also gives the best possible way to do so [1]. The Linear Quadratic Regulator is one such design methodology whereby quadratic performance indices involving the control signal and the state variables are minimized in an optimal fashion. Till date, PID controllers are used widely in industrial process control applications due to their simple structure, tuning, ease of applicability and reliability [2]. Some efforts have been directed towards tuning PID controller with LQR technique as in He *et al.* [3] and Yu and Hwang [4], considering the error and integral of error as the state variables. The LQR based technique has also been extended for tuning PID controllers for sluggish over-damped second order processes in [3] by canceling one of the real system poles with a zero of the PID controller. Thus, the approach, presented in [3] does not give the flexibility of tuning oscillatory processes by selecting the optimal controller gains via LQR for the three state variables i.e. error, its rate and integral, which is circumvented in this paper. The advantage of LQR based optimal PID controller tuning, over existing time domain integral performance index minimization based method [5] is also addressed in this work.

Fractional order controllers are gaining increasing interest in the research community due to recent hardware realizations, and scope for wide applicability and enhanced performance [6]. The FOPID or $PI^\lambda D^\mu$ controller, proposed by Podlubny [7] is an extension of conventional PID controllers and is gradually getting more importance in various process control applications. Due to its extra degrees of freedom the fractional order $PI^\lambda D^\mu$ controller has higher capability of enforcing control objectives than the conventional integer order PID controller. However the performance of such controller greatly depends on its tuning strategy. Several tuning methodologies have been proposed to tune $PI^\lambda D^\mu$ controllers like Ziegler-Nichols type [8], analytical rule based [9], stabilization based methods [10], time domain optimal tuning [11]-[13], frequency domain robust tuning [14]-[17], coupled time and frequency domain optimization based [18]-[19], optimization based dominant pole placement [20] etc. In this paper we propose a new LQR based formulation for the FOPID or $PI^\lambda D^\mu$ controllers for a class of fractional order processes with one non-integer order element ($\alpha$) as studied in [16]-[17]. The motivation of this work is to bridge the gap between the quadratic optimal control and fractional order control application using $PI^\lambda D^\mu$ structure to handle FO systems. LQR based FOPID controller design has been attempted in the present study, by formulating a non-commensurate order state space model while considering the error signal and its fractional order differ-integral as the state variables. In the present approach, the diagonal elements of the weighting matrix ($Q$) and weighting factor ($R$) are chosen as the decision variables of the optimization algorithm. The corresponding optimal state-feedback gains (as the $PI^\lambda D^\mu$ controller gains) for the three states are then found out by solving the Continuous Algebraic Riccati Equation (CARE). A performance index as the weighted summation of ITAE and

ISCO is subsequently used to judge the optimal time domain performance of the closed loop system and is minimized using a real coded genetic algorithm.

He et al. [3] first proposed the technique to find out the weighting matrix ($Q$) and weighting factor ($R$) from closed loop damping and frequency specifications. In this paper, we have adopted an approach to find out the optimum set of weighting matrices for the optimal regulator design. Indeed, optimal controller obtained with the LQR approach automatically minimizes the variation in state trajectories but does not always show acceptable closed loop time domain response and might often include high overshoot, oscillations etc. In order to achieve efficient tracking of the set point change, the weighting matrices should be chosen in such a manner that it meets some additional time domain optimality criteria in terms of overshoot, rise and settling time etc. Wang et al. [21] used a genetic algorithm based technique to find out the optimal set of co-efficients of the weighting matrices of LQR by minimizing another time domain optimality criterion. Poodeh et al. [22] also adopted a similar approach of finding weighting matrices with GA by the minimization of a custom cost function of steady-state error, maximum percentage of overshoot, rise time and settling time. The idea has been extended to FOPID controller design by Padula and Visioli [11] with IAE criterion and maximum sensitivity constraint and by Cao et al. [12] by weighted summation of IAE and ISCO. In the present study, the introduction of ITAE as the error index penalizes the chance of high oscillations at later stages and effectively reduces the rise time and settling time.

Like, the conventional integer order case, optimal control theory has also been extended for fractional order systems with fractional state variables by Agrawal [23]. Shafieezadeh et al. [24] have investigated the effect of fractional powers of the state variables along with the conventional optimal state feedback law. Tricaud and Chen [25], Agrawal [26], Biswas and Sen [27] formulated the fractional optimal control problem with a quadratic performance index involving the states and control law. Li and Chen [28], Tangpong and Agrawal [29], Biswas and Sen [30] also proposed similar quadratic performance index in matrix form for fractional optimal control problems that has been adopted in the present work. Saha et al. [31] studied LQR equivalence of dominant pole placement problem with FOPID controllers. Das et al. [32] studied the proposed the optimum weight selection based digital PID design for fractional order integral performance indices. This paper improvises over the available techniques by coupling the LQR theory with GA based time domain optimal $PI^\lambda D^\mu$ controller tuning, via an optimal choice of weighting matrix ($Q$) and weighting factor ($R$) of the CARE for an incommensurate fractional order state space formulation [6], while keeping the flexibility of choosing integro-differential orders of FOPID controller separately, unlike [10].

Rest of the paper is organized as follows. Section II discusses about the optimal state feedback approach for FOPID controller tuning and optimum selection of weighting matrices. Section III presents simulation study of the proposed controller with oscillatory and sluggish process. The paper ends with conclusions as section IV, followed by references.

## II. THEORETICAL FORMULATION FOR OPTIMAL FRACTIONAL ORDER CONTROLLER DESIGN

### A. State-feedback Approach for FOPID Controller Tuning for Fractional Order Systems

The classical LQR based optimal state-feedback controllers minimizes the quadratic cost function (1) with $\{x, u\}$ being the state variables and control actions respectively.

$$J = \int_0^\infty \left[ x^T(t)Qx(t) + u^T(t)Ru(t) \right] dt \tag{1}$$

Minimization of the integral performance index (1) leads to the solution of continuous time Riccati equation given by (2) to find out the state-feedback control law (3).

$$A^T P + PA - PBR^{-1}B^T P + Q = 0 \tag{2}$$
$$u(t) = -R^{-1}B^T Px(t) \tag{3}$$

In fact, the choice of weighting matrices ($Q$ and $R$) and also the integral and derivative orders of fractional states of the $PI^\lambda D^\mu$ controller ($\lambda$ and $\mu$) does not affect the optimal regulator formulation given by (1)-(3). For each choice of $\{Q, R\}$ it is possible to find out one optimal solution (controller gains) by solving the CARE given by (2). It is obvious that the closed loop performance changes a lot with variation in $\{Q, R, \lambda, \mu\}$. So, our objective is to find out the most optimal solution among all the optimal controllers that can by designed by LQR as in [21]-[22]. In order to do this, the optimal regulator formulation has been improved with another time domain optimality criterion, given by the weighted sum of ITAE and ISCO (discussed in the next subsection), since this criterion simultaneously minimizes the overshoot, increases the speed of time response and also reduces the control signal.

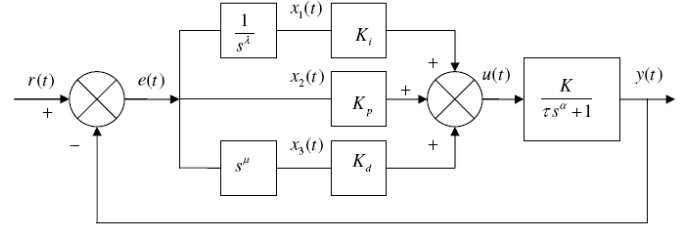

Figure 1. State feedback formulation of FOPID controller for FO-plants.

The formulation of the LQR based FOPID controller for controlling a class of fractional order plant [16]-[17] has been shown in Fig. 1, where the fractional differ-integrals ($\lambda$ and $\mu$) of the error signal has been considered as the state variables. The pseudo time-constant, dc gain and system order are denoted by $\{\tau, K, \alpha\}$. In [6], it is shown that the system shows highly oscillatory ($2 > \alpha > 1$) and highly sluggish ($1 > \alpha > 0$) open loop response, even with a simple first order like template due to the presence of fractional order dynamics of the plant.

In Fig. 1, if the system is excited with an external input $r(t)$ to get a control signal $u(t)$ and output $y(t)$, let us consider,

$$x_1(t) = D^{-\lambda}[e(t)], x_2(t) = e(t), x_3(t) = D^{\mu}[e(t)] \quad (4)$$

where, differ-integral operator $D = d/dt$ and $e(t)$ is the error signal. Therefore,

$$\frac{d^{\lambda}}{dt^{\lambda}}[x_1(t)] = x_2(t); \quad \frac{d^{\mu}}{dt^{\mu}}[x_2(t)] = x_3(t) \quad (5)$$

It has been reported in He et al. [3] that in the case of feedback control design, the external set-point does not affect the result and we can put $r = 0$. Thus the closed loop system is reduced to a regulator problem [1]. Clearly, for set-point $r = 0$, the error signal becomes $e = r - y = -y$. Thus, the output signal $y = -e = -x_2$. Now, from Fig. 1 it is clear that the fractional order process [16]-[17] is given by

$$\frac{Y(s)}{U(s)} = \frac{K}{\tau s^{\alpha} + 1} \Rightarrow -\frac{E(s)}{U(s)} = \frac{K}{\tau s^{\alpha} + 1} \quad (6)$$

$$\Rightarrow \tau D^{\alpha}[x_2(t)] + x_2(t) = -Ku(t)$$

Assuming zero initial condition for FO composition rule, the above equation reduces to

$$\tau D^{\alpha-\mu}x_3(t) + x_2(t) = -Ku(t) \quad (7)$$

$$\Rightarrow D^{\alpha-\mu}x_3(t) = -\frac{1}{\tau}x_2(t) - \frac{K}{\tau}u(t)$$

From (4), (5) and (7) we get,

$$\begin{bmatrix} \frac{d^{\lambda}}{dt^{\lambda}}x_1(t) \\ \frac{d^{\mu}}{dt^{\mu}}x_2(t) \\ \frac{d^{\alpha-\mu}}{dt^{\alpha-\mu}}x_3(t) \end{bmatrix} = \begin{bmatrix} 0 & 1 & 0 \\ 0 & 0 & 1 \\ 0 & -\frac{1}{\tau} & 0 \end{bmatrix}\begin{bmatrix} x_1(t) \\ x_2(t) \\ x_3(t) \end{bmatrix} + \begin{bmatrix} 0 \\ 0 \\ -\frac{K}{\tau} \end{bmatrix}u(t) \quad (8)$$

This is the representation of the above system in generalized incommensurate fractional order state-space template, i.e.

$$\frac{d^q x(t)}{dt^q} = Ax(t) + Bu(t) \quad (9)$$

Here, $\frac{d^q}{dt^q} = \begin{bmatrix} \frac{d^{\lambda}}{dt^{\lambda}} & \frac{d^{\mu}}{dt^{\mu}} & \frac{d^{\alpha-\mu}}{dt^{\alpha-\mu}} \end{bmatrix}^T$ and

$$x(t) = \begin{bmatrix} x_1(t) & x_2(t) & x_3(t) \end{bmatrix}^T = \begin{bmatrix} \frac{d^{-\lambda}e(t)}{dt^{-\lambda}} & e(t) & \frac{d^{\mu}e(t)}{dt^{\mu}} \end{bmatrix}^T$$

The system matrices are given by

$$A = \begin{bmatrix} 0 & 1 & 0 \\ 0 & 0 & 1 \\ 0 & -\frac{1}{\tau} & 0 \end{bmatrix}; B = \begin{bmatrix} 0 \\ 0 \\ -\frac{K}{\tau} \end{bmatrix} \quad (10)$$

Let us consider,

$Q = \begin{bmatrix} Q_1 & 0 & 0 \\ 0 & Q_2 & 0 \\ 0 & 0 & Q_3 \end{bmatrix}$ and $P = \begin{bmatrix} P_{11} & P_{12} & P_{13} \\ P_{12} & P_{22} & P_{23} \\ P_{13} & P_{23} & P_{33} \end{bmatrix}$ to solve the

CARE in (2). For a guess value of weighting matrix $Q$ and $R$, the six elements of the symmetric positive definite matrix i.e. $\{P_{11}, P_{22}, P_{23}, P_{12}, P_{13}, P_{23}\}$ can be solved out using MATLAB's Control System Toolbox [33] function *lqr()*. Therefore, the state feedback gain matrix ($F$) is obtained as:

$$F = R^{-1}B^T P = \frac{1}{R}\begin{bmatrix} 0 & 0 & -\frac{K}{\tau} \end{bmatrix}\begin{bmatrix} P_{11} & P_{12} & P_{13} \\ P_{12} & P_{22} & P_{23} \\ P_{13} & P_{23} & P_{33} \end{bmatrix}$$

$$= \frac{1}{R}\begin{bmatrix} \left(-\frac{K}{\tau}P_{13}\right) & \left(-\frac{K}{\tau}P_{23}\right) & \left(-\frac{K}{\tau}P_{33}\right) \end{bmatrix} \quad (11)$$

$$= \begin{bmatrix} -K_i & -K_p & -K_d \end{bmatrix}$$

Since the state variables are so chosen that it represents the error signal and its fractional differ-integrals, the design of the optimal state feedback regulator with it gives the $PI^{\lambda}D^{\mu}$ controller gains as the optimal state feedback gain matrix ($F$). The corresponding optimal control law is given by (12).

$$u(t) = -Fx(t) = -R^{-1}B^T Px(t)$$

$$= -\begin{bmatrix} -K_i & -K_p & -K_d \end{bmatrix}\begin{bmatrix} \frac{d^{-\lambda}e(t)}{dt^{-\lambda}} & e(t) & \frac{d^{\mu}e(t)}{dt^{\mu}} \end{bmatrix}^T \quad (12)$$

$$= K_p e(t) + K_i \underbrace{\int \cdots \int}_{\lambda-fold} e(t)dt + K_d \frac{d^{\mu}e(t)}{dt^{\mu}}$$

*B. Optimum Selection of Weighting Matrices*

The above LQR based formulation can be termed as optimal for a specific choice of the weighting matrices $Q$ and $R$. Indeed, the time domain performance is heavily affected for any arbitrary choice of the weighting matrices although the optimality is preserved. This is logical since the choice of weighting matrices determine the state feedback gains ($PI^{\lambda}D^{\mu}$ controller gains in this case) which directly affect the performance of the closed loop system. In order to handle this problem, a GA based stochastic optimization based technique is proposed by minimizing another time domain performance index $\tilde{J}$ (13) as studied in [5], [32] which tunes the elements of the weighting matrices i.e. $\{Q_1, Q_2, Q_3, R\}$ and integro-differential orders of the FOPID controller i.e. $\{\lambda, \mu\}$ while minimizing a weighted sum of ITAE and ISCO.

$$\tilde{J} = \int_0^{\infty}\left[w_1 \cdot t|e(t)| + w_2 \cdot u^2(t)\right]dt \quad (13)$$

Here, $w_1$ and $w_2$ are the corresponding weights of ITAE and ISCO and are considered to be same, so as to put equal impact on the error minimization criteria and small control signal. The rationale for using both these parameters in the objective function is to get a good time domain response and at the same time to limit the controller output to avoid actuator saturation and integral wind-up [2]. At a first glance this might seem as a redundant repetition since the LQR methodology already gives optimal values of the controller gains with the lowest cost. However, this is actually obtained for a specified value of the weighting matrices. When $Q$ and $R$ are varied, for each choice of weighting matrices the LQR would give an optimal gain with the lowest possible cost, but that does not necessarily imply a good time domain performance [21]-[22] with the LQR cost function (1). Also, for an optimal choice of weighting matrices ($Q$ and $R$) and differ-integral orders ($\lambda$ and $\mu$), the FOPID tuning problem becomes optimal due to the introduction of time domain optimal performance (13) as well as the classical optimal regulator (LQR) based approach (1), involving the fractional states.

The performance of the proposed approach is compared with the simple time domain optimal FOPID tuning technique of Cao *et al.* [12] and Das *et al.* [5] with the minimization of weighted ITAE and ISCO given as (13). This comparison verifies that alone the output optimality criteria (13) cannot ensure state optimality (1). Also, as far as low control signal is concerned, a $PI^\lambda D^\mu$ controller design with coupled output and state optimality criteria via LQR is a much powerful technique over the simple time-domain performance index minimization approaches as in [11]-[13].

### C. Optimization Algorithm Used for Controller Tuning

Genetic Algorithm is a stochastic optimization algorithm and has been widely employed in the tuning of PID controllers subject to the minimization of time domain cost function. Genetic algorithm has certain advantages over the classical gradient based optimization algorithms since they are stochastic in nature and are less susceptible to get trapped in the local minima within the search space. Initially a random population of genes (which is essentially a vector comprising of the decision variables) is chosen from the search space. They undergo reproduction, crossover and mutation to yield individuals with better fitness (lower $\tilde{J}$ value in this case). The individuals with higher fitness values have more probability of creating their copies in the next generation. This is termed as reproduction. Two parent individuals can do information interchange in a probabilistic fashion to create a child in the next generation. This process is known as crossover. In mutation a small part of the parent gene is randomly changed to yield a child. Another factor called the elite count is used which dictates the number of fittest individuals that would definitely go to the next generation. In this case, the population size is considered to be 20 and elite count as 2.

In the present simulation a crossover fraction of 0.8 and a mutation fraction of 0.2 are chosen which works well for a wide variety of problems [34]. Also similar to Zamani *et al.* [18], a high penalty is imposed when the choices of controller gains give an unstable closed loop response. The algorithm is terminated if the change in the objective function is lower than a specified tolerance level [34]. The variables that constitute the search space for the fractional order PID controller are $\{Q_1, Q_2, Q_3, R, \lambda, \mu\}$. The intervals of the search space for these variables are $\{Q_1, Q_2, Q_3, R\} \in [0, 100]$ and $\{\lambda, \mu\} \in [0, 2]$. In fact, with this search interval the number of state variables remains always three in the state space formulation (8) even though the integral and differential orders may take values higher than unity, representing faster time response and better closed stability respectively. The decisison variables are encoded as real values in the GA. The algorithm has also been run several times to ensure that the true global minima is found in the search space and only the best results with the lowest cost function (along with the corresponding decision variables) have been reported in the present study.

### III. ILLUSTRATIVE EXAMPLES

Two fractional order plants have been considered, as in (6), showing heavily oscillatory and sluggish open loop response. The two test plants have been chosen from [35], having different values of the system order as in [16]-[17]. These fractional order processes show sluggish and oscillatory open loop dynamics for $(\alpha < 1)$ and $(\alpha > 1)$ respectively.

### A. Simulation Study for an Oscillatory Fractional Order Process

The oscillatory system under consideration is represented by the following transfer function [35]:
$$G_1(s) = \frac{5}{1.11s^{1.7} + 1} \quad (14)$$
For the process (14), the LQR based tuning with the optimal selection of weighting matrices yields the $PI^\lambda D^\mu$ controller
$$C_1^{LQR}(s) = 0.726453 + \frac{0.692674}{s^{0.998773}} + 0.582319 s^{0.386624} \quad (15)$$
with the corresponding weighting matrices of performance index $J$ in (1) as
$$Q_1 = \begin{bmatrix} 0.474582 & 0 & 0 \\ 0 & 0.011476 & 0 \\ 0 & 0 & 0.01637 \end{bmatrix}; R_1 = 0.989131$$
The continuous Riccati solution in this case is given by
$$P_1 = \begin{bmatrix} 0.634755 & 0.398973 & 0.152102 \\ 0.398973 & 0.381525 & 0.15952 \\ 0.152102 & 0.15952 & 0.12787 \end{bmatrix}$$
Also, the optimal time domain tuning strategy proposed by Cao *et al.* [12] which only minimizes the cost function $\tilde{J}$ in (13) with genetic algorithm, yields another $PI^\lambda D^\mu$ controller as (16) for the oscillatory process (14)
$$C_1^{ITAE+ISCO}(s) = 0.100718 + \frac{0.93109}{s^{0.997477}} + 0.834496 s^{0.357018} \quad (16)$$

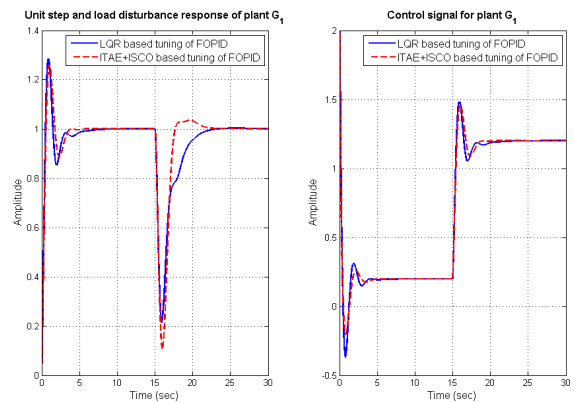

Figure 2. Response and control for plant $G_1$ with two $PI^\lambda D^\mu$ tuning strategies.

The simulated results of closed loop performances of plant (14) with the $PI^\lambda D^\mu$ controller (15) and (16), tuned with two different philosophies have been compared in Fig. 2. Also, Fig. 3 shows the state trajectories as the fractional differ-integrals of the error signal with the LQR based and conventional time domain error index and control signal based $PI^\lambda D^\mu$ controller tuning techniques.

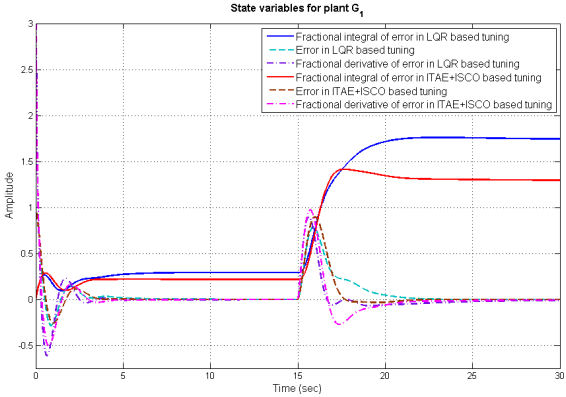

Figure 3. Time evolution of the state variables for tuning plant $G_1$.

### B. Simulation Study for a Sluggish Fractional Order Process

A sluggish process is considered next having the transfer function (17), similar to that in [35].

$$G_2(s) = \frac{5}{1.11s^{0.7}+1} \quad (17)$$

The LQR based tuning with GA based optimum weight selection by minimizing $\tilde{J}$ (13) yields the following $PI^\lambda D^\mu$ controller as (18)

$$C_2^{LQR}(s) = 1.900408 + \frac{2.302821}{s^{0.948591}} + 0.940017s^{0.017093} \quad (18)$$

with the corresponding weighting matrices as

$$Q_2 = \begin{bmatrix} 1.599235 & 0 & 0 \\ 0 & 0.012767 & 0 \\ 0 & 0 & 0.012018 \end{bmatrix}; R_2 = 0.301573$$

The continuous Riccati solution in this case is given by

$$P_2 = \begin{bmatrix} 1.458666 & 0.652811 & 0.154172 \\ 0.652811 & 0.441259 & 0.127231 \\ 0.154172 & 0.127231 & 0.062933 \end{bmatrix}$$

Also, the GA based minimization of (13) yields another optimal $PI^\lambda D^\mu$ controller as (19) for the sluggish process (17).

$$C_2^{ITAE+ISCO}(s) = 0.937303 + \frac{4.636422}{s^{0.949254}} + 0.030218s^{0.043881} \quad (19)$$

The closed loop performances and the control signals for plant (17) with controller (18) and (19) are shown in Fig. 4. Also, the corresponding error signal and its fractional differ-integrals have been shown as the state variables in Fig. 5.

It is evident from Fig. 2 and Fig. 4 that the load disturbance suppression is much better for the proposed LQR based tuning of FOPID controllers, than optimum ITAE and ISCO based tuning [12]. Also, it can be seen that the initial control signal is much lower for optimal weighting matrices based LQR tuning of FOPID controllers. Hence, the advantage of LQR based tuning is evident for $PI^\lambda D^\mu$ controllers also as reported with integer order PID controllers by He et al. [3]. Additionally the deviation in the state trajectories are enforced by the LQR technique which cannot be enforced by a simple weighted ITAE and ISCO based tuning like [12].

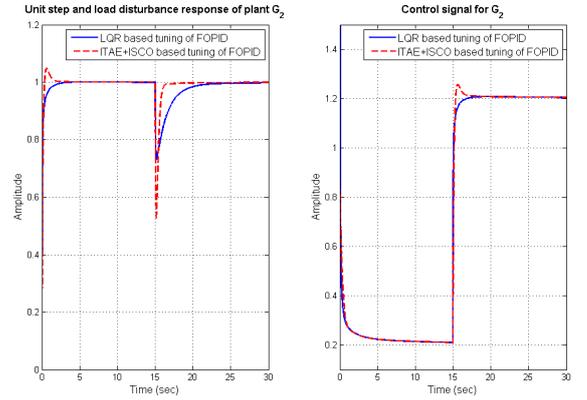

Figure 4. Response and control for plant $G_2$ with two $PI^\lambda D^\mu$ tuning strategies.

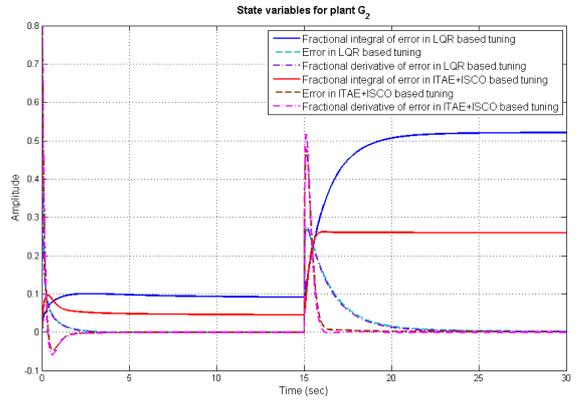

Figure 5. Time evolution of the state variables for tuning plant $G_2$.

### IV. CONCLUSION

LQR based improved fractional order $PI^\lambda D^\mu$ controller tuning has been proposed in this paper with optimal selection of weighting matrix $Q$ and weighting factor $R$. The optimal choice of the weighting matrices along with the fractional order integro-differential operators of the $PI^\lambda D^\mu$ controller have been obtained through real coded GA based minimization of a time domain performance indices, comprising of ITAE and ISCO. Thus, the proposed method preserves the state optimality of LQR and at the same time gives a low error index in the closed loop time response. This typical improvement enables the designer to obtain satisfactory closed loop response while also enjoying the benefits of LQR in the optimal $PI^\lambda D^\mu$ controller tuning. Simulation results show that the proposed techniques works well even for a special class of fractional order systems, with highly oscillatory and highly sluggish open loop response, in terms of low overshoot, better ability to suppress load disturbances and low initial controller effort for unit change in set-point. Future scope of work may include LQR based FOPID controller tuning for unstable and integrating fractional order systems.